# Optimization of Passive Chip Components Placement with Self-Alignment Effect for Advanced Surface Mounting Technology

Irandokht Parviziomran[a], Shun Cao[a], Haeyong Yang[b], Seungbae Park[c], Daehan Won[a]

[a]Department of System Science and Industrial Engineering, State University of New York, Binghamton, New York 139
[b]Koh Young Technology America,, Binghamton, New York 13902
[c]Department of Mechanical Engineering State University of New York, Binghamton, New York 13902

**Abstract**

Surface mount technology (SMT) is an enhanced method in electronic packaging in which electronic components are placed directly on soldered printing circuit board (PCB) and are permanently attached on PCB with the aim of reflow soldering process. During reflow process, once deposited solder pastes start melting, electronic components move in a direction that achieve their highest symmetry. This motion is known as self-alignment since can correct potential mounting misalignment. In this study, two noticeable machine learning algorithms, including support vector regression (SVR) and random forest regression (RFR) are proposed as a prediction technique to (1) diagnose the relation among component self-alignment, deposited solder paste status and placement machining parameters, (2) predict the final component position on PCB in $x$, $y$, and rotational directions before entering in the reflow process. Based on the prediction result, a non-linear optimization model (NLP) is developed to optimize placement parameters at initial stage. Resultantly, RFR outperforms in terms of prediction model fitness and error. The optimization model is run for 6 samples in which the minimum Euclidean distance from component position after reflow process from ideal position (i.e., the center of pads) is outlined as 25.57 ($\mu m$) regarding defined boundaries in model.

*Keywords:* Passive chip Component self alignment; optimal placement parameters setting; machine learning; predictive modeling; optimization-based modeling; pick and placement process; surface mount technology

## 1. Introduction

Surface mount technology (SMT) is an advanced method in electronic packaging in which solder paste provides electrical connection between surface mount components (SMCs) and electronic circuit. Three main operations are accomplished to produce assembled electronic circuit, including stencil printing process (SPP), pick and placement process (P&P), and reflow soldering. At first, solder paste is deposited onto the surface of printing circuit board (PCB). Then SMCs are placed automatically with the aim of P&P machine. Finally, solder joints which are formed by employed heat during reflow process, attach SMCs on PCB permanently. SMT enhances electronic packages by utilization of smaller packages, smaller circuit boards [1], smaller components, etc. On one side, miniaturization makes component placement more challenging since more accurate component placement is required with smaller size of packages with larger lead counts [1]. On the other side, mounted components move during reflow soldering process unintentionally because of the forces that acting on component when solder paste starts wetting. These phenomena is known as the 'self-alignment' and originated from the surface tension of molted solder that empowers component to move in a direction that achieves its equilibrium state [2]. So, to obtain the knowledge of how and how much chip components move during reflow concerning their initial status such as components geometry (miniaturized chip components), designed pad geometry, deposited solder paste status, and placement strategy, would be critical and brings the opportunity of optimizing placement parameter setting that reduces components positional offsets after reflow soldering.

However, numerous studies have investigated chip component self-alignment capability from theoretical [1], simulation and numerical [1, 3] models' standpoints, there is no generalized data-driven model in literature to address practical challenges of self-alignment [4]. Lv et al. provided a comprehensive survey on machine learning application in SMT [4]. Based on this survey, there is no research to employ applied machine learning method in self-alignment [4] while applied machine learning methods privilege over conventional statistical methods in SMT [5]. Despite a recently published paper by Marktinek et al. that used a one-layer neural network to predict component position after reflow regarding of its placement offsets in $x$, $y$, and rotational directions [6] while neglecting other factors. Moreover, the overarching goal of developing a self-alignment prediction model is preventing potential defects such as tombstoning, overhanging, etc. before entering reflow process. Particularly, tombstoning happens when SMCs are lifted from a pad on one end in which electrical connection between SMCs and PCB would be broken. Overhanging refers to the situation when SMCs land off from their pad in $x$ and $y$ directions. Accurate SMCs placement regarding deposited solder paste volume and offsets can reduce the possibility of tombstoning and overhanging. So, an optimization model that optimizes placement setting could be accomplished with the prediction model to avoid such defects. For the best of our knowledge, no optimization model in the literature addresses placement machining setting in P&P process. Hence, this enhances the necessities of: (1) obtaining insight of factors that contribute in component self-alignment, (2) developing a generalized data-driven model that addresses different type and size of chip components, (3) outlining machine learning algorithm that predicts chip components position after reflow soldering, (4) determining the best placement machining parameters that enhances the quality of finished electronic package. This study considers 6 types of passive chip components (3 resistors and 3 capacitors). A generalized optimization-prediction model could be a breakthrough in surface mount assembly (SMA) line since on one side, all previous studies considered only one or two types of passive chip components self-alignment [3, 5, 6] and on the other side, there is no optimization-prediction model that addresses placement machining strategy.

From the extensive research and based on the domain knowledge, 13 variables are selected to train the prediction models, and three targets are defined as components offsets after reflow process, naming as post offset $x$, $y$, and rotation ($\theta$) (i.e., post refers to component position after reflow process). Two remarkable machine learning algorithms including support vector regression (SVR) and random forest regression (RFR) are used. Then, we build a non-linear optimization (NLP) model to determine the optimized component placement setting by setting the controllable parameters in the P&P machine. The NLP model is solved with evolutionary strategy (ES) because of its complexity.

The rest of this paper is organized as follows: Section 2 presents the literature related to components self-alignment; the built prediction models, optimization model and ES are discussed in Section 3 followed by results in Section 4; the conclusions and future work of this research are provided in Section 5.

## 2. Literature review

Chip components capability of being self-aligned has been investigated in several studies with the means of dynamic fluid concepts which describe the scientific reason of this motion. Restoring forces, mainly force originating from surface tension, acts on components once solder paste is in its fluid state. This force drags chip components in a direction that achieves its symmetry. The primary challenge regarding this phenomenon is obtaining insight into factors that contribute to components movement. For this matter, the experimental result is accomplished with force models as well as simulation models in the literature. A dynamic force model proposed by [1] demonstrates the effect of pad geometry, chip metallization and dimensions, solder volume, and placement offsets on components motion during the reflow process. This study considered one type of capacitor named as 1206. Based on the results, smaller pad lengths, smaller pad gaps, larger solder volume, and smaller metallization provided better self-alignment for this capacitor [1]. The effect of different types of solder, different types of reflow technology, different solder paste volume and different placement angles on self-alignment have been studied in [7] in which the importance of solder volume on components self-alignment along with utilization of lead-free solders have been outlined [7]. Moreover, Liukkonen et al. compared component offsets before and after reflowing with leaded and lead-free processes [2]. However, they proposed that lead-free process has more variation in self-alignment comparing leaded process [2], other potential factors like paste volume and offsets are neglected from their study. Chip mass, die tilt, and solder volume variations and distributions have been considered in the regression model by [8] to optimize static equilibrium conditions of a flip-chip. They concluded that the solder volume variation has more significant influence on the chip standoff height (i.e., in the $z$ direction) than the chip lateral alignment accuracy (i.e., in $x$ and $y$ directions) [8].

## 3. Methodology

In this section, a comprehensive description of the experiment and collected data are presented and contributing factors under this study are introduced. Furthermore, brief descriptions of applied machine learning techniques and their justification for this study are also presented. Finally, the optimization model is discussed along with the metaheuristic algorithm that is considered to solve the proposed optimization problem.

## 3.1. Data description

The experiment is designed to assemble 6 passive chip component types, including 3 resistors and 3 capacitors in 3 size categories as R1005 and C1005 (1mm×0.5mm), R0603 and C0603 (0.6mm×0.3mm), and R0402 and C0402 (0.4mm×0.2mm). 660 placements is considered for each component type (6 different components) which means 3940 placements in total. All components are placed horizontally on two corresponding pads (For detail, see Fig. 1). Longer and shorter sides of the component are considered as parallel with $x$ and $y$ directions respectively, and rotation is defined as a circular movement on $x$-axis. The red dot in the Fig. 1 indicates the center of two pads which is considered as a reference point (i.e., {0, 0}). At first stage, the solder paste is deposited on the PCB and various factors including paste volume ratio on pad 1 and 2, and paste offsets ($x$ and $y$) on pad 1 and 2 are measured with the means of advanced solder paste inspection (SPI) machine. Fig. 1(a) presents solder paste positional factors in $xy$-plane. Volume average ratio and volume difference ratio are defined as the average and difference of deposited volume on pad 1 and pad 2 respectively. Then, components are placed with intentional offsets in $x$, $y$, and rotational ($\theta$) directions with the aid of

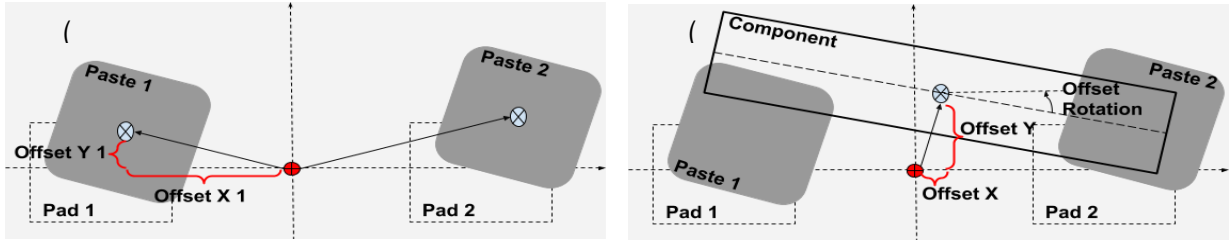

Fig. 1. positional relation of pads, solder pastes and component in $xy$-plane. (a) paste offsets description; (b) placement offset description

P&P machine (see Fig. 1(b) for description of placement offsets). Finally, components are permanently attached on PCB by nitrogen reflow oven with given thermal profile following the justification of used lead-free solder paste guidelines. Table 1 and Table 2 present considered categorical and continuous factors with a brief statistical description from this experiment respectively.

Table 1. brief explanation and levels of selected categorical variables

| Categorical Variables | Description | Name | Level 1 | Level 2 | Level 3 | Level 4 | Level 5 |
|---|---|---|---|---|---|---|---|
| Component Sizes ($1000\mu m^2$) | $Component\ length \times Component\ width$ | $\lambda_1$ | 80 | 180 | 500 | N/A | N/A |
| Component Types | Resistor / Capacitor | $\lambda_2$ | R | C | N/A | N/A | N/A |
| Pad Sizes ($1000\mu m^2$) | $Pad\ length \times Pad\ width$ | $\lambda_3$ | 44 | 102 | 280 | N/A | N/A |
| Pad Gaps ($\mu m$) | Gap btw. two pads | $\lambda_4$ | 160 | 260 | 250 | 450 | 460 |

Table 2. brief explanation and levels of selected continuous variables

| Continuous Variables | Description | Name | Max | Min | $(\mu, \sigma)$ |
|---|---|---|---|---|---|
| Volume Avg. (%) | Avg. of printed volume ratio on 2 corresponding pads | $\gamma_1$ | 154.77 | 46.32 | (95.25, 16.62) |
| Volume Diff. (%) | Diff. of printed volume ratio on 2 corresponding pads | $\gamma_2$ | 96.40 | -91.48 | (2.21, 27.85) |
| Paste Offset X 1 ($\mu m$) | Center of deposited solder paste on pad 1 in x direction | $\gamma_3$ | 698.81 | 188.54 | (403.96, 159.05) |
| Paste Offset Y 1 ($\mu m$) | Center of deposited solder paste on pad 1 in y direction | $\gamma_4$ | 316.63 | 18.86 | (129.22, 64.60) |
| Paste Offset X 2 ($\mu m$) | Center of deposited solder paste on pad 2 in x direction | $\gamma_5$ | -76.96 | -412.58 | (-216.29, 90.52) |
| Paste Offset Y 2 ($\mu m$) | Center of deposited solder paste on pad 2 in y direction | $\gamma_6$ | 316.63 | 8.13 | (130.00, 62.87) |
| Pre Offset X ($\mu m$) | Center of component in x direction after placement | $\chi_1$ | 316.91 | -37.15 | (123.16, 78.56) |
| Pre Offset Y ($\mu m$) | Center of component in y direction after placement | $\chi_2$ | 264.57 | -97.88 | (61.38, 58.78) |
| Pre Offset Rotation ($deg.$) | Component circular movement on x-axis | $\chi_3$ | 24.78 | -32.90 | (-0.12, 2.98) |

## 3.2. Support vector regression (SVR)

The SVR is a prediction technique that is able to solve complex nonlinear regression problems by mapping input features into high-dimensional space $F$, wherein they are linearly correlated with the problem target. For a given training set with sample size $n$; $\{(x_1, y_1), \ldots, (x_n, y_n)\} \subset X \times \mathbb{R}$ where $X$ denotes the space of input features ($X = \mathbb{R}^d$, $d$ is number of features) and $\mathbb{R}$ is the set of real numbers for continuous target $y$, the SVR method fits hyperplane $F(X) = \langle w, x \rangle + b$, in which most data points fall on this plane [9]. Note that $\langle \cdot, \cdot \rangle$ denotes the dot product of any two arbitrary vectors. The distance of $\varepsilon$ and $-\varepsilon$ from defined hyperplane is defined as marginal tolerance range to penalize prediction point out of this boundary with a predefined penalty ratio of $C\sum_{i=1}^{n}(\xi_i + \xi_i^*)$ in which $\xi_i$ and $\xi_i^*$ denotes the amount of deviation from defined hyperplane for $i = 1, \ldots, n$, where $n$ is data size. The SVR estimation function and corresponding constraints are shown in Eq. (1) and (2) in which the objective is minimizing the norm of hyperplane weight ($w$) and penalty to balance the hyperplane and its tolerance margins.

$$\min \quad \frac{1}{2}\|w\|^2 + C\sum_{i=1}^{n}(\xi_i + \xi_i^*) \tag{1}$$

$$s.t. \quad \begin{cases} y_i - \langle w, x_i \rangle - b \leq \varepsilon + \xi_i \\ -y_i + \langle w, x_i \rangle + b \leq \varepsilon + \xi_i^* \\ \xi_i, \xi_i^* \geq 0 \;\; for\; i = 1, \ldots, n \end{cases} \tag{2}$$

In this study, the ε-insensitive loss function is considered along with the linear kernel. For the SVR setting, we set $\varepsilon = 0.1$ as the distance of boundaries from hyperplane and $C = 1$ as the penalty ratio.

### 3.3. Random forest regression (RFR)

The RFR is an ensemble prediction technique which is built based on a collection of randomized regression trees. In RFR, each randomly selected tree employs random split on feature space with the objective of reducing prediction error. Finally, RFR, acquires the mean of individual outputs to get more stable prediction [10]. For a given training set with sample size $n$; $\{(x_1, y_1), \ldots, (x_n, y_n)\} \subset X \times \mathbb{R}$, where $X$ denotes the space of input features ($X = \mathbb{R}^d$, $d$ is the number of features) and $\mathbb{R}$ is the set of real numbers for continuous target $y$, a random forest is defined as a hierarchical tree-structured predictor [10] in which: $F(x) = \frac{1}{J}\sum_{j=1}^{J} f_j(x)$, where $J$ is the number of trees in the forest and $f_j(x)$ is the estimation function of $j$th tree that is trained on a random split of input features [10]. For a decision tree with $M$ splitting nodes, the feature space would be splitting into $M$ regions as $R_m$ in which: $f_j(x) = \sum_{m=1}^{M} b_m \varphi(x, R_m)$; $for\; j = 1, \ldots, J$, where $b_m$ is the corresponding constant for each region. $\varphi(x, R_m)$ is the binary decision function that shows whether input feature $x$ is selected (i.e., $\varphi(x, R_m)=1$ if $x \in R_m$) or not ($\varphi(x, R_m) = 0$ if $x \notin R_m$). In this study, the number of trees is initialized as 50 and trees are fully grown.

### 3.4. Passive chip component placement NLP model

The placement parameters that are needed to be optimized are the last 3 continuous variables ($\chi_d$; $d \in D = \{1,2,3\}$) shown in Table 2, and a combination of them represents a placement setting for specific component placement in mounting process. The first 6 continuous variables ($\gamma_i$; $i \in I = \{1, \ldots, 6\}$) in Table 2, are solder paste properties corresponding each pad. Finally, the 4 categorical variables ($\lambda_j$; $j \in J$; $J = \{1, \ldots, 4\}$) in

Table 1 indicates a specific component directory and its corresponding pad directory that is designed to be placed on PCB. The solder paste properties-related variables and component and pad directory-related variables are considered as features in training SVR and RFR prediction models to capture the effect of different components and solder paste status on components self-alignment. These variables are given variables in optimization model. The main objective of the NLP model is to reduce the Euclidean distance of component position after reflow soldering from reference point ($R_x, R_y$) (see Eq. (3)). Rotational offset ($\theta$) is also addressed as a hard constraint in Eq. (4). The proposed NLP model considers optimization for component position in $x$, $y$, and $\theta$ directions with a user predefined threshold for each (see Eq. (4)-(6)). The reason behind this is that no placement would accrue with zero offset concerning ideal position in real production line. The list of notations used to drive the objective function and constraints of the proposed NLP optimization model is shown in Table 3. The functional values are estimated based on the trained SVR and RFR prediction models, as shown:

Table 3. list of notations

| | Notation | Name | Description |
|---|---|---|---|
| Indices and Sets | $d \in D$ | | Set of decision variables in optimization model; placement parameters setting in $x$, $y$, and rotational directions |
| | $i \in I$ | | Set of given continuous variables in prediction model; solder paste properties-related factors |
| | $j \in J$ | | Set of given categorical variables in prediction model; component and pad directory-related factors |
| | $k \in K$ | | Set of level of categorical variables; $L = \{1, 2, 3, 4, 5\}$ |
| Variables | $\chi$ | | Decision variable; placement parameters setting |
| | $\gamma$ | | Given continuous variable in prediction model |
| | $\lambda$ | | Given categorical variable |
| Parameters | $(x, y, \theta)$ | | Component positional offsets in $x$, $y$, and rotational directions |
| | $\tau$ | | Threshold of prediction-optimization constraint |
| | $L_d, U_d$ | | Lower and upper bounds for placement parameter $d$, respectievely |
| Function | $F$ | | Prediction function (SVR and RFR) |

$$\min \quad \sqrt{\| R_x - F_x(\lambda_{jk}, \gamma_i, \chi_d) \|^2 + \| R_y - F_y(\lambda_{jk}, \gamma_i, \chi_d) \|^2} \qquad (3)$$

$$s.t. \quad |F_\theta(\lambda_{jk}, \gamma_i, \chi_d)| \leq \tau_\theta \qquad (4)$$

$$|F_x(\lambda_{jk}, \gamma_i, \chi_d)| \leq \tau_x \qquad (5)$$

$$|F_y(\lambda_{jk}, \gamma_i, \chi_d)| \leq \tau_y \qquad (6)$$

$$L_d \leq \chi_d \leq U_d; d \in D; i \in I; j \in J; k \in K \qquad (7)$$

The thresholds (i.e., $\tau_\theta, \tau_x, \tau_y$ ) in Eq. (4)-(6) are considered as customer defined expectation limits and are defined as $\tau_\theta = 2$ degree, $\tau_x = 20\%$ of pad length $(\mu m)$, $\tau_y = 20\%$ of pad width $(\mu m)$ in this study, based on experts' opinion. Furthermore, since the aim of this model is to drive a placement strategy for a given component and pad directory and solder paste properties, decision variables are limited with lower and upper boundaries in which lower bound presents center of two pads and upper bound presents center of two printed pastes for decision variables $\chi_1$ and $\chi_2$ such that $(L_{\chi_1}, U_{\chi_1}) = (R_x, (\gamma_3 + \gamma_5)/2)$ and $(L_{\chi_2}, U_{\chi_2}) = (R_y, (\gamma_4 + \gamma_6)/2)$. Concerning upper bound of $\chi_3$, the slope of two deposited solder pastes is calculated as $(L_{\chi_3}, U_{\chi_3}) = (R_\theta, arctan((\gamma_4 - \gamma_6)/(\gamma_3 - \gamma_5)))$.

*3.5. Solution approach*

Since the proposed optimization model involves probabilistic SVR and RFR prediction functions in its formulation (non-linear functions), it is hard to solve with deterministic optimization. For this matter, a population-based metaheuristic method, evolutionary strategy (ES), is used [11] to find the optimal placement setting for proposed NLP problem. Basically, ES saves the best from a population to amend next generation [11]. So, in the end, the best population survives with ES. ES is chosen in this research rather than other population-based metaheuristics because the optimization problem of this research consists of continuous and categorical variables with real-valued decision variables. So, ES makes the configuration of proposed NLP simpler in terms of encoding real-valued variables [11]. Moreover, ES converges to an optimal or near optimal solution faster than other methods such as genetic algorithm (GA) [11]. The formulation in [11] is modified for proposed NLP as shown in Algorithm 1. The proposed algorithm employs normal distribution $N(0, \sigma)$ to produce offspring and $(\mu, \lambda)$-ES is used as a selection strategy in which best children are selected to create the next generation [11].

Algorithm 1. ES to retrieve optimal placement parameters; $\mu = 5, \lambda = 10, \sigma = 0.5, G = 10$

---

Initialize $P^\mu$ as the number of parents in population, $\lambda$ as the number of offspring, $G$ as the number of generations, $N(0, \sigma)$ as mutation operator. Next generation population is selected based on $(\mu, \lambda)$-ES.
Randomly selected population $P$ of size $\mu$
**while** $g = \{1, ..., G\}$ **do**
    **for** $m = 1, ..., \mu$ **do**
        **select** parent $p$ with a uniform probability from $P^\mu$
        **for** $o = 1, ..., \lambda$ **do**
            **mutate** parent p by $\sigma$ to form $\lambda$ offspring
            **Check** $bounderies == Ture$ & $constraint == True$ in proposed NLP
        **end**
    **end**
    **return** $P^{\mu \times \lambda}$ and let $I$ be a solution in $P^{\mu \times \lambda}$
    **calculate** the probability of successful mutation by $\frac{1}{5}$ rule and update $\sigma$
    **for** $t = 1, ..., \mu \times \lambda$ **do**
        **evaluate** fitness function in the proposed NLP
    **end**
    **return** $f^{\mu \times \lambda}$ objective values
    from $f^{\mu \times \lambda}$ **select** top $\mu$ best fitness and corresponding $S^\mu$ solutions from $P^{\mu \times \lambda}$
    **replace** $P^\mu == S^\mu$
    $g += 1$
**end**
**return** best solution from the final population $P^\mu$

## 4. Results and discussion

The result of experiment is divided in 70: 10: 20 as training: validating: testing sets in order to: (1) train SVR and RFR prediction models on the training set and tune the parameters of models with validating set; (2) test SVR and RFR models with the unseen testing set. Moreover, the coefficient of determination ($R^2$) is employed to measure the fitness of learning algorithm on training set along with root-mean-squared error (RMSE) which is used to evaluate the performance of the prediction model on testing set. Table 4 presents the $R^2$ value of training set and an RMSE value of testing set for SVR and RFR prediction models.

Table 4. Testing accuracy measures of SVR and RFR prediction models

| Prediction Model | Post Offset $x$ ($\mu m$) | | Post Offset $y$ ($\mu m$) | | Post Offset $\theta$ ($deg.$) | |
|---|---|---|---|---|---|---|
| | $RMSE(\mu m)$ | $R^2$ | $RMSE(\mu m)$ | $R^2$ | $RMSE(\mu m)$ | $R^2$ |
| SVR | 18.32 | 0.38 | 16.65 | 0.42 | 1.61 | 0.02 |
| RFR | 15.48 | 0.94 | 12.63 | 0.95 | 1.56 | 0.87 |

Based on the prediction results, RFR outperforms in terms of both model fitness on the training set and prediction error on the testing set. The model is fitted with $R^2$ value of 94%, 95%, and 87% for components offset after reflow soldering in $x$, $y$, and $\theta$ directions respectively. Moreover, RFR model predicts components position after reflow soldering with 15.48 ($\mu m$), 12.65 ($\mu m$), and 1.56 (deg.) prediction error in $x$, $y$, and $\theta$ directions respectively. At the next step, the probabilistic function of RFR is employed in the optimization model. Basically, RFR's mathematical mapping function from input feature space to output space, is used in optimization model. So, for any given solution in feature space such as solder paste average volume ratio, the optimization model searches for optimal component position before reflow soldering ($x, y, \theta$) in which it has minimum offsets after reflow soldering (see NLP formulation in section 3.4). To show the result, we run the optimization model for six samples, each component type once. In Fig. 2, the pink triangle shows actual result from experiment and violet rectangle presents optimized result for component offset before and after reflow soldering respectively. It is worthy to mention that components position is not necessary equal before and after reflow soldering because components move during reflow process by the self-alignment effect. As mentioned before, the aim of this study is finding the best placement in respect of such motion during reflow process. Based on the results, the optimization model is terminated with the objective value of 25.57 (i.e., component distance from ideal position after reflow process) for all samples regarding proposed boundaries and thresholds in section 3.4. As shown in Fig. 2, optimal placement parameter (i.e., violet rectangle) varies dramatically from actual placements in our experiment. To achieve the least offsets after reflow soldering, it is suggested to place components with initial offsets of (0, 50) ($\mu m$) in $x$ and $y$ directions and (-2, 2) (deg.) for rotation. Moreover, rerunning the experiment with optimized placement parameters could verify this optimization model. However, errors originating from the prediction model and SMT machinery operations would affect the final accuracy regardless of the optimized placement.

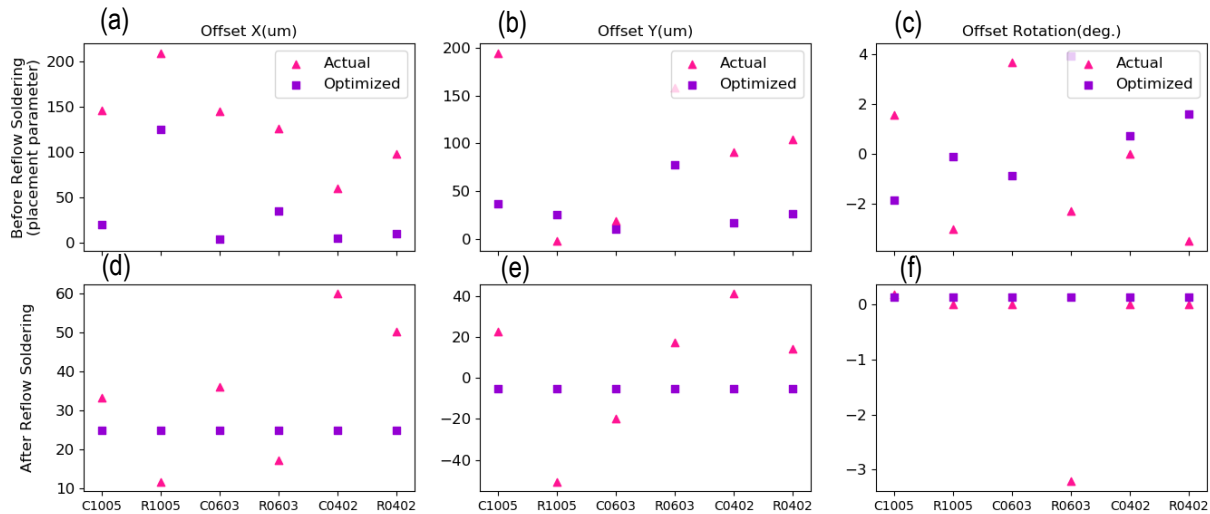

Fig. 2. experiment result ( actual) vs. optimized components offsets (a) before reflow soldering in $x$ direction; (b) before reflow soldering in $y$ direction; (c) before reflow soldering in rotational direction; (d) after reflow soldering in x direction; (e) after reflow soldering in y direction; (f) after reflow soldering in rotational direction

## 5. Conclusion and future work

This study employed 13 potential factors that contribute in components movement during reflow soldering and applied two magnitude machine learning approaches as SVR and RFR to predict passive chip components position after reflow soldering in $x$, $y$ and rotational directions. Then an NLP model is developed to optimize component placement setting including placement in $x$, $y$ and rotational directions with the objective of minimizing the probability function of proposed prediction models. Based on the result, RFR outperforms in terms of training fitness and prediction error and is chosen for the optimization model. The optimization results, which are examined on 6 samples, show the objective value as 25.57 ($\mu m$) that indicates the minimum Euclidean distance from component position after reflow process from ideal position (i.e., center of pads) regarding considered boundaries and thresholds in the model. Moreover, placement parameters for each sample is also retrieved. Based on the result for these 6 samples, optimal placement setting is between (0, 50) ($\mu m$) in $x$ and $y$ directions and between (-2, 2) (deg.) for rotation based on specific given factors including solder paste average and difference volume ratio and solder paste offsets ($x$ and $y$). It is also recommended to verify optimization results with real experiment.

However, the proposed prediction-optimization model is the first in component self-alignment, but four main areas would be recommended as future works. Specifically, it is to: (1) develop the multi-target prediction model that considers all targets simultaneously; (2) reduce prediction error with the aid of advanced machine learning techniques such as feature selection, hyperparameter tuning, etc.; (3) consider all prediction targets in objective function of optimization model; (4) design stochastic version of the optimization model instead of deterministic model to address errors origination from prediction model as well as random error in the machinary operations.